%% file: AnotherApproachToGetDerivativeOfOddPower.tex
\documentclass[12pt,letterpaper,oneside,reqno]{amsart}
\usepackage{amsfonts}
\usepackage{amsmath}
\usepackage{amssymb}
\usepackage{amsthm}
\usepackage{float}
\usepackage{mathrsfs}
\usepackage{colonequals}
\usepackage[font=small,labelfont=bf]{caption}
\usepackage[left=1in,right=1in,bottom=1in,top=1in]{geometry}
\usepackage[pdfpagelabels,hyperindex,colorlinks=true,linkcolor=blue,urlcolor=magenta,citecolor=green]{hyperref}
\usepackage{setspace}
\usepackage{array}
\usepackage{etoolbox}
\apptocmd{\sloppy}{\hbadness 10000\relax}{}{}
\newcommand \coeffA [3][A] {{\mathbf{#1}} \sb{#2,#3}}

\newtheorem{thm}{Theorem}[section]
\newtheorem{example}[thm]{Example}
\newtheorem{definition}[thm]{Definition}
\newtheorem{notation}[thm]{Notation}

\title[Another approach to get derivative of odd-power]
{Another approach to get derivative of odd-power}
\author[Petro Kolosov]{Petro Kolosov}
\address{Software Developer, DevOps Engineer}
\email{kolosovp94@gmail.com}
\urladdr{https://kolosovpetro.github.io}
\date{\today}
\hypersetup{
    pdftitle={Another approach to get derivative of odd-power},
    pdfsubject={
        Derivative, Partial derivatives, Partial differential equations, Polynomials
    },
    pdfauthor={Petro Kolosov},
    pdfkeywords={
        Derivative, Partial derivatives, Partial differential equations, Polynomials
    }
}
\begin{document}
    \begin{abstract}
        \input{sections/01_abstract}
    \end{abstract}

    \maketitle

    \tableofcontents

    \section{Introduction and Main Results} \label{sec:introduction}
    \input{sections/02_introduction_and_main_results}
    \input{sections/examples/01_example}
    \input{sections/examples/02_example}
    \input{sections/examples/03_example}

    \section{Conclusions}\label{sec:conclusions}
    \input{sections/03_conclusions}

    \section{Verification of the results}\label{sec:verification-of-the-results}
    \input{sections/04_verification}

    \bibliographystyle{unsrt}
    \bibliography{AnotherApproachToGetDerivativeOfOddPower}
    \noindent \textbf{Version:} \input{sections/version}

\end{document}

%% file: sections/01_abstract.tex
In this manuscript, we provide and discuss another approach to get derivative of odd-power
such that is based on an identity in partial derivatives in terms of polynomial function $f_y$ defined as
\[
    f_{y} (x, z) = \sum_{k=1}^{z} \sum_{r=0}^{y} \mathbf{A}_{y,r} k^r (x-k)^r
\]
where $x, z\in \mathbb{R}$, $y$ is fixed constant $y \in \mathbb{N}$ and $\mathbf{A}_{y,r}$ are real coefficients.

%% file: sections/02_introduction_and_main_results.tex
This manuscript provides another approach to get derivative of odd-power,
that is an approach based on partial derivatives of the polynomial function $f_y(x,z)$ defined as
\[
    f_{y} (x, z) = \sum_{k=1}^{z} \sum_{r=0}^{y} \mathbf{A}_{y,r} k^r (x-k)^r
\]
where $x, z\in \mathbb{R}$, $y$ is fixed constant $y \in \mathbb{N}$ and $\mathbf{A}_{y,r}$ are real coefficients.
The essence of the approach we discuss is build
on an identity in terms of sum of partial derivatives of the polynomial function $f_{y}$.
The function $f_{y}$ is defined by the main results of the manuscript~\cite{kolosov_2022}
that explains an odd-power in a form as follows
\begin{equation}
    n^{2m+1} = \sum_{k=1}^{n} \sum_{r=0}^{m} \coeffA{m}{r} k^r (n-k)^r
    \label{eq:odd-power-identity}
\end{equation}
where $m$ is fixed constant $m\in\mathbb{N}$, $n \in \mathbb{N}$ and $\coeffA{m}{r}$ are real coefficients defined
recursively, see~\cite{kolosov2016link}.
We define the function $f_{y}$ such that based on the identity~\eqref{eq:odd-power-identity}
with the only difference that values of $n, m$ in the right part of~\eqref{eq:odd-power-identity}
appear to be parameters of the function $f_{y}$.
In contrast to the equation~\eqref{eq:odd-power-identity}, upper bound $n$ of the sum $\sum_{k=1}^{n}$ turned into fixed
function's parameter $y$ as well.
Let the function $f_{y}$ be defined as follows
\begin{definition} (Polynomial function $f_{y}$.)
    \begin{equation}
        f_{y} (x, z) = \sum_{k=1}^{z} \sum_{r=0}^{y} \coeffA{y}{r} k^r (x-k)^r
        \label{eq:definition-f}
    \end{equation}
\end{definition}
where $x, z\in \mathbb{R}$ and $y$ is constant $y \in \mathbb{N}$.
Note that for every $x\in\mathbb{R}$ and constant $y\in\mathbb{N}$ the polynomial identity satisfies
\begin{equation*}
    f_{y} (x, x) = x^{2y+1}
\end{equation*}
At first glance, the definition~\eqref{eq:definition-f} might look complex, so in order to clarify
the function $f_y$ and polynomials it produces, let there be a few examples.
Substituting the values of $y=1,2,3$ to the function $f_y$ we get the following polynomials in $x, z$
\begin{align*}
    f_{1} (x, z) &= 3 x z - 3 z^2 + 3 x z^2 - 2 z^3 \\
    f_{2} (x, z) &= 5 x^2 z - 15 x z^2 + 15 x^2 z^2 + 10 z^3 - 30 x z^3 + 10 x^2 z^3 +
    15 z^4 - 15 x z^4 + 6 z^5 \\
    f_{3} (x, z) &= -7 x z + 14 x^2 z + 7 z^2 - 42 x z^2 + 35 x^3 z^2 + 28 z^3 - 140 x^2 z^3 + 70 x^3 z^3 + 175 x z^4 \\
    &- 210 x^2 z^4 + 35 x^3 z^4 - 70 z^5 + 210 x z^5 - 84 x^2 z^5 - 70 z^6 + 70 x z^6 - 20 z^7
\end{align*}
These polynomials are obtained by rearranging the sums in the definition~\eqref{eq:definition-f} as
\[
    f_{y} (x, z) = \sum_{r=0}^{y} \coeffA{y}{r} \left[ \sum_{k=1}^{z} k^r (x-k)^r \right]
\]
So that part $\sum_{k=1}^{z} k^r (x-k)^r$ is polynomial
in $x, z$ calculated using Faulhaber's formula~\cite{beardon1996sums}.
According to the main topic of the current manuscript, it provides another approach to get derivative of odd-power.
Therefore, we define odd-power function we work in the context of.
The odd-power function $g_y$ is a function defined as follows
\begin{definition}(Odd-power function $g_y$.)
    \begin{equation*}
        g_{y}(x) = x^{2y + 1}
    \end{equation*}
\end{definition}
where $x\in \mathbb{R}$ and $y$ is constant $y\in \mathbb{N}$.
The Interesting part is that odd-power function $g_{y} (x)$ may be obtained
as a partial case of the function $f_y$ for $z=x$.
Also, the ordinary derivative of odd-power $\frac{d}{dx} g_{y}$ evaluate in point $u\in\mathbb{R}$
may be obtained as a sum of partial derivatives of $f_y$ evaluate in point $(u,u)$.
We explain this further in the manuscript.
One more important thing remains to conclude is to define partial derivative's notation.
More precisely, the following notation for partial derivatives
is used across the manuscript and remains unchanged
\begin{notation} (Partial derivative.)
    Let be a function $\beta (x_1, x_2, \dots, x_n)$ defined over the real space $\mathbb{R}^n$.
    We denote partial derivative of the function $\beta$ with respect to $x_i$ as follows
    \begin{equation*}
        \beta^{'}_{x_i}
        = \lim_{\Delta x_i \to 0}
        \frac{\beta (x_1, x_2, \dots, x_i + \Delta x_i, \dots, x_n) - \beta (x_1, x_2, \dots, x_n)}{\Delta x_i}
    \end{equation*}
\end{notation}
Partial derivative of the function $\beta_{x_i}$ with respect to $x_i$
evaluate in point $(y_1, y_2, \dots, y_n) \in \mathbb{R}^n$ is denoted as follows
\begin{equation*}
    \beta^{'}_{x_i} (y_1, y_2, \dots, y_n)
\end{equation*}
Moreover, partial derivative $\beta^{'}_{x_i}$ evaluate in point $(y_1, y_2, \dots, y_n) \in \mathbb{R}^n$ plus
partial derivative $\beta^{'}_{x_j}$ evaluate in point $(y_1, y_2, \dots, y_n) \in \mathbb{R}^n$
is equivalent to the sum of partial derivatives $\beta^{'}_{x_i} + \beta^{'}_{x_j}$
evaluate in point $(y_1, y_2, \dots, y_n) \in \mathbb{R}^n$ and to be denoted as
\begin{equation*}
    \beta^{'}_{x_i} (y_1, y_2, \dots, y_n)
    + \beta^{'}_{x_j} (y_1, y_2, \dots, y_n)
    = [\beta^{'}_{x_i} + \beta^{'}_{x_j}] (y_1, y_2, \dots, y_n)
\end{equation*}
So that now we can switch our focus back to the functions $g_{y}$ and $f_{y}$.
Therefore, the following theorem in terms of partial derivatives
reflects the relation between the ordinary derivative of odd-power
function $g_{y}$ and function $f_{y}$
\begin{thm}
    \label{thm:main-theorem}
    Let be a fixed point $v\in \mathbb{N}$, then ordinary derivative $\frac{d}{dx} g_v (u)$ of the odd-power function $g_v(x) = x^{2v + 1}$
    evaluate in point $u\in\mathbb{R}$ equals to partial derivative $(f_{v})^{'}_{x} (u, u)$ evaluate in point $(u, u)$ plus
    partial derivative $(f_{v})^{'}_{z} (u, u)$ evaluate in point $(u, u)$
    \begin{equation}
        \frac{d}{dx} g_v (u) = (f_{v})^{'}_{x} (u, u) + (f_{v})^{'}_{z} (u, u)
        \label{eq:odd-exponential-identity}
    \end{equation}
\end{thm}
In particular, it follows that for every pair $u \in \mathbb{R}, v \in \mathbb{N}$ an identity holds
\begin{align*}
(2v+1)
    u^{2v} &= (f_{v})^{'}_{x} (u, u) + (f_{v})^{'}_{z} (u, u) \\
    &= [(f_{v})^{'}_{x} + (f_{v})^{'}_{z}](u,u)
\end{align*}
that is also an ordinary derivative of odd-power function $t^{2v+1}, \; v\in \mathbb{N}, \; v=\mathrm{const}$
evaluate in point $u\in\mathbb{R}$,
therefore
\begin{align*}
    \frac{d}{dt} t^{2v+1} (u) &= (f_{v})^{'}_{x} (u, u) + (f_{v})^{'}_{z} (u, u) \\
    &= [(f_{v})^{'}_{x} + (f_{v})^{'}_{z}](u,u)
\end{align*}
To summarize and clarify all about the theorem~\ref{thm:main-theorem}, we provide a few examples that show
an application of it.

%% file: sections/examples/01_example.tex
\begin{example}
    \normalfont
    Theorem~\ref{thm:main-theorem} example for $x\in\mathbb{R}, \; z\in \mathbb{R}$ and $y=1$.
    Consider the explicit form of the function $f_{1} (x, z)$, that is
    \[
        f_1(x, z) = 3 x z - 3 z^2 + 3 x z^2 - 2 z^3
    \]
    Therefore, the partial derivative $(f_1)^{'}_{x}$ with respect to $x$ equals to
    \[
        (f_1)^{'}_{x} = \lim_{d \to 0} \frac{3 d z + 3 d z^2}{d} = 3 z + 3 z^2
    \]
    Consider the partial derivative $(f_1) ^{'}_{z}$ with respect to $z$, that is
    \begin{align*}
    (f_1) ^{'}_{z}
        &= \lim_{d \to 0} \left[\frac{-3 d^2 - 2 d^3 + 3 d x + 3 d^2 x - 6 d z - 6 d^2 z + 6 d x z - 6 d z^2}{d} \right] \\
        &= \lim_{d \to 0} \left[ -3 d - 2 d^2 + 3 x + 3 d x - 6 z - 6 d z + 6 x z - 6 z^2 \right] \\
        &=3 x - 6 z + 6 x z - 6 z^2
    \end{align*}
    Summing up both partial derivatives $(f_1)^{'}_{x}$ and $(f_1)^{'}_{z}$, we get
    \begin{align*}
    (f_1)
        ^{'}_{x} + (f_1)^{'}_{z} = 3 x - 3 z + 6 x z - 3 z^2
    \end{align*}
    Evaluating in point $(u, u)$ yields
    \begin{align*}
        \frac{d}{dt} t^{3} (u) = [(f_1)^{'}_{x} + (f_1)^{'}_{z}] (u,u)  = 3 u^2
    \end{align*}
    That confirms the results of the theorem~\ref{thm:main-theorem}.
\end{example}

%% file: sections/examples/02_example.tex
\begin{example}
    \normalfont
    Theorem~\ref{thm:main-theorem} example for $x\in\mathbb{R}, \; z\in \mathbb{R}$ and $y=2$.
    Consider the explicit form of the function $f_{2} (x, z)$, that is
    \[
        f_2 (x, z) = 5 x^2 z - 15 x z^2 + 15 x^2 z^2 + 10 z^3 - 30 x z^3 + 10 x^2 z^3 + 15 z^4 - 15 x z^4 + 6 z^5
    \]
    Therefore, the partial derivative $(f_2) ^{'}_{x}$ with respect to $x$ equals to
    \begin{align*}
    (f_2) ^{'}_{x} &= \lim_{d \to 0} \left[ 5 d z + 10 x z - 15 z^2 + 15 d z^2 + 30 x z^2 - 30 z^3 + 10 d z^3 +
        20 x z^3 - 15 z^4 \right] \\
        &= 10 x z - 15 z^2 + 30 x z^2 - 30 z^3 + 20 x z^3 - 15 z^4
    \end{align*}
    Consider the partial derivative $(f_2)^{'}_{z}$ with respect to $z$, that is
    \begin{align*}
    (f_2)
        ^{'}_{z}
        &= 5 x^2 - 30 x z + 30 x^2 z + 30 z^2 - 90 x z^2 + 30 x^2 z^2 + 60 z^3 - 60 x z^3 + 30 z^4
    \end{align*}
    Summing up both partial derivatives $(f_2)^{'}_{x}$ and $(f_2)^{'}_{z}$, we get
    \begin{align*}
    (f_2)
        ^{'}_{x} + (f_2)^{'}_{z} &= 5 x^2 - 20 x z + 30 x^2 z + 15 z^2 - 60 x z^2 + 30 x^2 z^2 + 30 z^3 - 40 x z^3 + 15 z^4
    \end{align*}
    Evaluate in point $(u, u)$ yields
    \begin{equation*}
        \frac{d}{dt} t^{5} (u) = [(f_2)^{'}_{x} + (f_2)^{'}_{z}] (u,u) = 5 u^4
    \end{equation*}
    That confirms the results of the theorem~\ref{thm:main-theorem}.
\end{example}

%% file: sections/examples/03_example.tex
\begin{example}
    \normalfont
    Theorem~\ref{thm:main-theorem} example for $x\in\mathbb{R}, \; z\in \mathbb{R}$ and $y=3$.
    Consider the explicit form of the function $f_{3} (x, z)$, that is
    \begin{align*}
        f_3 (x, z) &= -7 x z + 14 x^2 z + 7 z^2 - 42 x z^2 + 35 x^3 z^2 + 28 z^3 -140 x^2 z^3 + 70 x^3 z^3 + 175 x z^4 \\
        &- 210 x^2 z^4 + 35 x^3 z^4 -70 z^5 + 210 x z^5 - 84 x^2 z^5 - 70 z^6 + 70 x z^6 - 20 z^7
    \end{align*}
    Therefore, the partial derivative of $(f_3) ^{'}_{x}$ with respect to $x$ equals to
    \begin{align*}
    (f_3) ^{'}_{x} &= -7 z + 28 x z - 42 z^2 + 105 x^2 z^2 - 280 x z^3 + 210 x^2 z^3 + 175 z^4 - 420 x z^4 \\
        &+ 105 x^2 z^4 + 210 z^5 - 168 x z^5 + 70 z^6
    \end{align*}
    Consider the partial derivative $(f_3) ^{'}_{z}$ with respect to $z$, that is
    \begin{align*}
    (f_3) ^{'}_{z} &= -7 x + 14 x^2 + 14 z - 84 x z + 70 x^3 z + 84 z^2 - 420 x^2 z^2 + 210 x^3 z^2 + 700 x z^3 \\
        &- 840 x^2 z^3 + 140 x^3 z^3 - 350 z^4 + 1050 x z^4 - 420 x^2 z^4 - 420 z^5 + 420 x z^5 - 140 z^6
    \end{align*}
    Summing up both partial derivatives $(f_3)^{'}_{x} (x, z)$ and $(f_3)^{'}_{z} (x, z)$, we get
    \begin{align*}
    (f_3) ^{'}_{x} + (f_3)^{'}_{z} &= -7 x + 14 x^2 + 7 z - 56 x z + 70 x^3 z + 42 z^2 - 315 x^2 z^2 + 210 x^3 z^2 \\
        &+ 420 x z^3 - 630 x^2 z^3 + 140 x^3 z^3 - 175 z^4 + 630 x z^4 - 315 x^2 z^4 - 210 z^5 \\
        &+ 252 x z^5 - 70 z^6
    \end{align*}
    Evaluate in point $(u,u)$ yields
    \begin{equation*}
        \frac{d}{dt} t^{3} (u) = [(f_3) ^{'}_{x} + (f_3)^{'}_{z}] (u,u) = 7 u^6
    \end{equation*}
    That confirms the results of the theorem~\ref{thm:main-theorem}.
\end{example}

%% file: sections/03_conclusions.tex
In this manuscript, we have reviewed an approach to get ordinary derivative of odd-power
using an identity in partial derivatives of the function $f_y$ evaluate in fixed point $(u,u) \in \mathbb{R}^2.$,
that is described by the theorem~\ref{thm:main-theorem}.
The main results of the manuscript can be validated
using Mathematica programs available online at~\cite{kolosov2022another}.

%% file: sections/04_verification.tex
As it is stated in conclusions, there is a possibility to validate the main results of this manuscript using
Wolfram Mathematica.
Therefore, a complete guide to validate the main results and formulae is attached as well.
Mathematica package source file is available online under the folder \texttt{mathematica}, see~\cite{kolosov2022another}.
The following expressions could be verified:
\begin{itemize}
    \item The function $f_y(x,z)$ for any constant argument $y\in\mathbb{N}$ using mathematica method \texttt{f[x, y, z]} e.g.
    \[
        \texttt{f[x, 1, z]} = 3 x z - 3 z^2 + 3 x z^2 - 2 z^3
    \]
    \item Partial derivative $(f_y)^{'}_x$ for any constant argument $y\in\mathbb{N}$ using mathematica
    method \texttt{DerivativeFByX[x, y, z]}
    \[
        \texttt{DerivativeFByX[x, 1, z]} = 3 z + 3 z^2
    \]
    \item Partial derivative $(f_y)^{'}_z$ for any constant argument $y\in\mathbb{N}$ using mathematica
    method \texttt{DerivativeFByZ[x, y, z]}
    \[
        \texttt{DerivativeFByZ[x, 1, z]} = 3 x - 6 z + 6 x z - 6 z^2
    \]
    \item Theorem~\ref{thm:main-theorem} for any constant argument $y\in\mathbb{N}$
    \begin{align*}
        \texttt{DerivativeFByX[x, 1, z]} + \texttt{DerivativeFByZ[x, 1, z]} &= 3 x - 3 z + 6 x z - 3 z^2 \\
        \texttt{DerivativeFByX[u, 1, u]} + \texttt{DerivativeFByZ[u, 1, u]} &= 3 u^2
    \end{align*}
\end{itemize}

%% file: sections/version.tex
\texttt{0.2.4-tags-v0-2-3.1+tags/v0.2.3.8e0e8a3}